\documentclass{amsart}
\usepackage{amssymb, amsthm, amscd, amsmath}

\newtheorem{theorem}{Theorem}[section]

\newtheorem{prop}[theorem]{Proposition}

\theoremstyle{definition}

\newtheorem{example}[theorem]{Example}

\theoremstyle{remark}
\newtheorem{remark}[theorem]{Remark}

\begin{document}

\def\CM{Cohen-Macaulay }
\def\deg{\mathop{\mathrm{deg}}\,}
\def\Proj{\mathop{\mathrm{Proj}}\,}
\def\Spec{\mathop{\mathrm{Spec}}\,}
\def\Cl{\mathop{\mathrm{Cl}}\,}
\def\Pic{\mathop{\mathrm{Pic}}\,}
\def\HH{\underline{\mathrm{Hom}}}
\def\H{\underline{H}}
\def\HE{\underline{\mathrm{Ext}}}
\def\HT{\underline{\otimes}}
\def\HR{\mathop{\mathrm{Hom}_R}}
\def\HK{\mathop{\mathrm{Hom}_K}}
\def\ER{\mathop{\mathrm{Ext}_R}}
\def\sym{\mathcal{R}_{\it{s}}}

\title{Cyclic covers of rings with rational singularities}

\author{Anurag K. Singh}

\address{Department of Mathematics, University of Utah, Salt Lake City,
155~S. 1400~E., UT~84112-0090} 

\curraddr{Mathematical Sciences Research Institute, 1000 Centennial Drive,
\#5070, Berkeley, CA~94720-5070}

\email{asingh@msri.org}

\thanks{This manuscript is based on work supported in part by the National
Science Foundation under Grant No. DMS 0070268. I would like to thank the
referee for a careful reading of the manuscript, and for helpful suggestions.}

\copyrightinfo{2002}{American Mathematical Society}


\dedicatory{Dedicated to Professor Phillip Griffith on the occasion of his
sixtieth birthday}

\begin{abstract} We examine some recent work of Phillip Griffith on \'etale
covers and fibered products from the point of view of tight closure theory.
While it is known that cyclic covers of Gorenstein rings with rational
singularities are Cohen-Macaulay, we show this is not true in general in the
absence of the Gorenstein hypothesis. Specifically, we show that the canonical
cover of a $\mathbb Q$-Gorenstein ring with rational singularities need not be
Cohen-Macaulay. \end{abstract}

\maketitle
\section{Introduction}

The theory of cyclic covers has found frequent use in commutative algebra in
recent years. In \cite{griffith-torsion} Griffith and Weston applied this
theory in their study of torsion elements in the kernel of certain divisor
class group homomorphisms. More recently, there has been a substantial amount
of interest in cyclic covers arising from tight closure theory and its
connections with the study of singularities in characteristic zero, see
\cite{sm-vanish}, in particular, \S\,4. While these connections have been
illuminated and greatly advanced by studying the notions of F-rationality and
F-regularity that arise from tight closure theory, the property of F-regularity
remains something of a mystery. It is known that a Gorenstein ring is F-regular
precisely if it is F-rational, and it is in the Gorenstein setting that the
property of F-regularity is perhaps best understood. A frequent theme has been
to extend results that are known for Gorenstein rings to the class of $\mathbb
Q$-Gorenstein rings, and it is here that cyclic covers make their appearance.
As an example, Hochster and Huneke showed that a Gorenstein integral domain of
characteristic $p>0$ is a direct summand of every module-finite extension ring
precisely if it is F-regular, \cite{HHjalg}; the question arises whether this
is true without the Gorenstein hypothesis, and the (albeit partial) answer to
this now available uses cyclic covers, see \cite{splinter}. In a geometric
vein, Hara's theorem that a ring of characteristic zero with log terminal
singularities has F-regular type, uses a cyclic cover argument to reduce to the
Gorenstein case where log terminal singularities agree with rational
singularities, see \cite{Hara}. It should be mentioned that rings of dimension
two with rational singularities are $\mathbb Q$-Gorenstein by the work of
Lipman, \cite{lipman}, and so it is only natural to consider their cyclic
canonical covers, which have the advantage of being Gorenstein rings. 

In \cite{wat-dim2} Watanabe showed that a cyclic cover of an F-regular ring is
also F-regular and, as a consequence, that it is Cohen-Macaulay. One of the
goals of this paper is to show that under the weaker hypothesis that a ring $R$
is F-rational (or that it has rational singularities), a cyclic cover of $R$
need not be a Cohen-Macaulay ring, even in the most natural setting where the
cyclic cover is precisely the canonical cover of a $\mathbb Q$-Gorenstein ring
$R$. Our main result is:

\begin{theorem} For all $d \ge 3$ there exists an $\mathbb N$-graded $\mathbb
Q$-Gorenstein ring $R$ of dimension $d$, with an isolated rational singularity
at its homogeneous maximal ideal, such that the canonical cover of $R$ is not a
\CM ring. \end{theorem}

In section \ref{phil} we provide a tight closure approach to some results from
Griffith's recent work \cite{griffith}. The question under consideration is the
following: let $f: X \to Y$ be a finite \'etale morphism of smooth projective
varieties over an algebraically closed field of characteristic zero. If $Y$ is
arithmetically Cohen-Macaulay, does this property ascend to $X$? More
precisely, if $\oplus_{n \ge 0} H^0(Y, {\mathcal L}^{\otimes n})$ is a
Cohen-Macaulay homogeneous coordinate ring for $Y$, where $\mathcal L$ is an
ample line bundle on $Y$, is the corresponding homogeneous coordinate ring 
$\oplus_{n \ge 0} H^0(X, f^* {\mathcal L}^{\otimes n})$ of $X$ also
Cohen-Macaulay? Similarly, one may ask the following question: let $(R,m)
\hookrightarrow (S,n)$ be a module-finite extension of local or $\mathbb
N$-graded rings which is \'etale on the punctured spectrum, $\Spec R - \{ m
\}$. If $R$ is Cohen-Macaulay, is the ring $S$ also Cohen-Macaulay? The answer
is negative as established by examples of Griffith, and we show how a very
general and interesting source of such examples comes from tight closure
theory.

\section{Graded rings and local cohomology}\label{graded}

We shall follow the notation and conventions of \cite{GW}. For the convenience
of the reader we record some standard facts about graded local cohomology
modules and Segre products.

By an ${\mathbb N}$-{\it graded ring} $R$, we will always mean a ring
$R=\oplus_{n \ge 0} R_n$, finitely generated over a field $R_0=K$. We use $m_R$
or $m$ to denote the homogeneous maximal ideal $R_+$ of $R$. By a homogeneous
system of parameters for $R$ we mean a sequence of homogeneous elements of $R$
whose images form a system of parameters for $R_m$. In specific examples
involving homomorphic images of polynomial rings, lowercase letters will denote
the images of the corresponding variables, the variables being denoted by
uppercase letters.

For a graded $R$-module $M$, the $i$\,th graded component of $M$ will be
denoted by $[M]_i$. We use $M(n)$ to denoted the graded $R$-module with
$[M(n)]_i = [M]_{n+i}$ for all $i \in {\mathbb Z}$. For graded $R$-modules $M$
and $N$, we may consider the graded $R$-module $\HH_R(M,N)$ whose $i$\,th
graded component $[\HH_R(M,N)]_i$ is the abelian group consisting of all {\it
graded\/} (i.e., degree preserving) R-linear homomorphisms from $M$ to $N(i)$. 
It is an easy verification that 
$$ 
\HH_R(M,N(n)) = \HH_R(M,N)(n) = \HH_R(M(-n),N).
$$ 
Graded $R$-modules form an abelian category with enough injectives. For an
integer $i \ge 0$, the functor $\HE_R^i(M,-)$ denotes the $i$\,th derived 
functor of $\HH_R(M,-)$. If $M$ is a finitely generated graded $R$-module, then 
$$ 
\HH_R(M,N) = \HR(M,N) \ \text{ and } \ \HE^i_R(M,N) = {\ER}^i(M,N) 
$$ 
as underlying $R$-modules, where $N$ is any graded $R$-module. For a
homogeneous ideal $I$ of $R$ and a graded $R$-module $M$, the graded local
cohomology modules of $M$ with support in $I$ (more precisely, with support in 
the closed set $V(I) \subseteq \Spec R$) are defined as 
$$ 
\H^i_I(M) = \varinjlim_n \, \HE^i_R (R/I^n, M). 
$$
Since the modules $R/I^n$ are finitely generated, $\H^i_I(M)$ and $H^i_I(M)$
agree as underlying $R$-modules. If $\dim R = d$ and $x_1, \dots, x_d$ is a 
homogeneous system of parameters for $R$, the module $\H_m^d(R)$ may be 
identified with 
$$ 
\varinjlim_t R /(x_1^t,\dots, x_d^t)R
$$ 
where the maps are induced by multiplication by the element $x_1 \cdots x_d \in
R$. Under the grading on $\H_m^d(R)$, we have
$$
\deg [z + (x_1^t,\dots, x_d^t)] = \deg z - t \sum_{i=1}^d \deg x_i
$$ 
where $z$ is a homogeneous element of $R$. Goto and Watanabe defined $a(R)$,
the {\em $a$-invariant} of $R$, as the highest integer $a$ such that
$[\H_m^d(R)]_a$ is nonzero. 

The injective hull of $K$ in the category of graded $R$-modules will be denoted
by $\underline E_R(K)$ and this is also the injective hull of $K$ as an
$R$-module. For a graded $R$-module $M$, its graded $K$-dual is 
$M^* = \HH_K(M,K)$ where
$$
[M^*]_i = \HK(M,K(i)) = \HK(M(-i),K).
$$
With this notation it turns out that $\underline{E}_R(K) = R^*$, and so
$$
\HH_R(M, \underline{E}_R(K)) \cong \HH_K(M,K) = M^*.
$$
The {\it graded canonical module}\/ of $R$ is 
$$
\omega_R = (\H_m^d(R))^* \cong \HH_R(\H_m^d(R), \underline{E}_R(K)).
$$ 
If the ring $R$ is Gorenstein, it is easy to see that there is a graded
isomorphism $\omega_R \cong R(a)$ where $a$ is the $a$-invariant of the ring 
$R$.

\subsection*{Segre products.}

Let $A$ and $B$ be $\mathbb N$-graded rings over a field $A_0 = B_0 = K$. The
{\it Segre product} of $A$ and $B$ is the ring
$$
R = A \# B = \oplus_{n \ge 0} A_n \otimes_K B_n.
$$
This has a natural grading in which $R_n = A_n \otimes_K B_n$. If $M$ and $N$ 
are graded modules over $A$ and $B$ respectively, then their Segre product is 
the graded $R$-module
$$
M \# N = \oplus_{n \in {\mathbb Z}} [M]_n \otimes_K [N]_n
\quad \text{ with } \quad [M \# N]_n = [M]_n \otimes_K [N]_n. 
$$
We assume next that $A$ and $B$ are normal rings over an algebraically closed
field $K$. In this case $R = A\#B$ is also a normal ring since it is a direct 
summand of the normal ring $A\otimes_K B$. For reflexive modules $M$ and $N$ 
over $A$ and $B$ respectively, we have the K\"unneth formula for local 
cohomology, \cite[Theorem 4.1.5]{GW}:
\begin{align*}
\H^k_{m_R}(M\# N) \ \cong & \left( M\#\H^k_{m_B}(N) \right) \ \oplus \ 
\left( \H^k_{m_A}(M)\#N \right) \\
& \oplus \bigoplus_{i+j=k+1} \left( \H^i_{m_A}(M) \# \H^j_{m_B}(N) \right)
\quad \text{ for all }\quad k \ge 0.
\end{align*}
If $\dim A = r \ge 1$ and $\dim B = s \ge 1$, the above formula shows that 
$\dim R = r+s-1$. Furthermore, if $r\ge 2$ and $s\ge2$, we get 
$$
\H^{r+s-1}_{m_R}(R) \ \cong \ \H^r_{m_A}(A) \# \H^s_{m_B}(B), 
$$
and applying $\HH_K(-,K)$, this shows that $\omega_R\cong \omega_A \# \omega_B$.

In the paper \cite{chow}, Chow established necessary and sufficient conditions
under which the Segre product of two graded \CM rings is Cohen-Macaulay. Let $U
\hookrightarrow {\mathbb P}^r$ and $V \hookrightarrow {\mathbb P}^s$ be
projective varieties with homogeneous coordinate rings $A$ and $B$
respectively. The Segre product $A \# B$ is a homogeneous coordinate ring for
the Segre embedding $U \times V \hookrightarrow {\mathbb P}^{rs+r+s}$, and the
question whether $U \times V \hookrightarrow {\mathbb P}^{rs+r+s}$ is
projectively \CM is precisely the question whether $A \# B$ is a \CM ring.
Chow's work was motivated by a question of Dwork and Ireland that arose from
studying Zeta functions of algebraic varieties. In \cite[\S 14]{HRinv} Hochster
and Roberts showed that under mild hypotheses on $A$ and $B$, Chow's results
can be obtained by an application of the K\"unneth formula for sheaf
cohomology. Also, Goto and Watanabe established such a result using the
K\"unneth formula for local cohomology. Their result is:

\begin{theorem}\cite[Theorem 4.2.3]{GW} 
Let $A$ and $B$ be $\mathbb N$-graded \CM rings of dimension at least two. If
$a(A) < 0$ and $a(B) < 0$, then the Segre product $A \# B$ is \CM and 
$a(A \# B) < 0$. Conversely if $A \# B$ is \CM and the rings $A$ and $B$ 
contain nonzero homogeneous elements of all positive degrees, then $a(A) < 0$ 
and $a(B) < 0$. 
\label{ainv-seg} \end{theorem}

\subsection*{Rational coefficient Weil divisors.} We recall some notation and
results from \cite{De, wat-dem, wat-dim2}.

By a {\it rational coefficient Weil divisor}\/ (or a $\mathbb Q$-{\it divisor})
on a normal projective variety $X$, we mean a $\mathbb Q$-linear combination of
codimension one irreducible subvarieties of $X$. If $D = \sum n_iV_i$ where
$n_i \in {\mathbb Q}$ and $V_i$ are distinct irreducible subvarieties, we set
$[D]= \sum [n_i]V_i$, where $[n]$ denotes the greatest integer less than or
equal to $n$, and define ${\mathcal{O}}_X(D) ={\mathcal{O}}_X([D])$.

Let $D=\sum(p_i/q_i)V_i$ where the integers $p_i$ and $q_i$ are relatively
prime, $q_i > 0$, and the subvarieties $V_i$ are distinct. The {\it fractional
part\/} of $D$ is defined as $D' = \sum((q_i-1)/q_i)V_i$. The reason for this
definition is that we then have $-[-nD] = [nD+D']$ for any integer $n$, and
this will be of use when we use Serre duality with $\mathbb Q$-divisors.

For an ample ${\mathbb Q}$-divisor $D$ (i.e., such that $nD$ is an ample 
Cartier divisor for some $n \in {\mathbb N}$), the {\it generalized section 
ring} corresponding to $D$ is
$$
R=R(X,D)=\oplus_{n\ge 0} H^0(X,{\mathcal {O}}_X(nD)).
$$
With this notation, Demazure's result (\cite[3.5]{De}) states that every 
${\mathbb N}$-graded normal ring $R$ is isomorphic to a generalized section 
ring $R(X,D)$ for some ample $\mathbb Q$-divisor $D$ on $X = \Proj R$.

Let $X$ be a smooth projective variety of dimension $d$ with canonical divisor
$K_X$, and let $D$ be an ample $\mathbb Q$-divisor on $X$. Let $\omega_R$ 
denotes the graded canonical module of the ring $R=R(X,D)$ and $\omega_R^{(i)}$ 
its $i$\,th symbolic power. In the papers \cite{wat-dem} and \cite{wat-dim2}
Watanabe showed that
\begin{align*} 
[\H_m^i(R)]_n & \cong H^{i-1}(X,{\mathcal O}_X(nD)) 
\quad \text{ for all } \quad i \ge 2, \\
[\omega_R^{(i)}]_n & \cong H^0(X,{\mathcal O}_X(i(K_X+D')+nD)), \quad 
\text{ and } \\
[\H^{d+1}_m(\omega_R^{(i)})]_n & \cong H^d(X,{\mathcal O}_X(i(K_X+D')+nD)).
\end{align*}

\subsection*{Rational singularities.} 

We shall say that a normal ring $R$ essentially of finite type over a field of
characteristic zero has {\it rational singularities}\/ if the affine scheme
$\Spec R$ has only rational singularities. We recall a theorem of Watanabe:

\begin{theorem}\cite[Theorem 2.2]{wat-kstar} Let $R$ be a normal $\mathbb
N$-graded ring which is finitely generated over a field $R_0 = K$ of
characteristic zero. Then $R$ has rational singularities if and only if the 
following two conditions are satisfied:
\begin{enumerate}
\item the open set $\Spec R - \{m_R\}$ has only rational singularities, and 

\item the ring $R$ is Cohen-Macaulay with $a(R) < 0$.
\end{enumerate}\label{ratsing-ainv} \end{theorem}

Note that condition $(1)$ above is satisfied if the ring $R$ has an isolated
singularity at its homogeneous maximal ideal $m_R$. An application of this
theorem, which we will have use for later, is the following:

\begin{prop} Let $D$ be an effective $\mathbb Q$-divisor on projective space
${\mathbb P}^d$ over a field of characteristic zero, and consider the
generalized section ring $R=R({\mathbb P}^d,D)$. If $\Spec R - \{m_R\}$ has at
most rational singularities, then $R$ has rational singularities.
\label{effective}\end{prop}

\begin{proof} By construction $R$ is a normal ring of dimension $d+1$, and 
since 
$$
[\H_m^i(R)]_n \cong H^{i-1}({\mathbb P}^d,{\mathcal O}_{{\mathbb P}^d}(nD))=0
$$
for all $2 \le i \le d$, for all $n \in \mathbb Z$, the ring $R$ must be
Cohen-Macaulay. It remains only to verify that $a(R) < 0$. Using $^*$ to denote
the Serre dual, 
$$
[\H_m^{d+1}(R)]_n \cong H^d({\mathbb P}^d,{\mathcal O}_{{\mathbb P}^d}(nD))
\cong H^0({\mathbb P}^d,{\mathcal O}_{{\mathbb P}^d}(K_{{\mathbb P}^d}+D'-nD))^*
$$
which is zero for all $n \ge 0$ since $D$ is effective.
\end{proof}

\section{Tight closure}

For the theory of tight closure, developed by Hochster and Huneke, we refer the
reader to the papers \cite{HHstrong, HHjams, HHbasec} and \cite{HHjalg}. In our
applications of this theory, we will usually be restricting ourselves to the
class of $\mathbb N$-graded rings, and we summarize some results which will
suffice for our needs. The reader can find more general forms of several of
these results in the papers mentioned above.

\begin{theorem} Let $R$ be an $\mathbb N$-graded ring over a perfect field 
$R_0 = K$ of characteristic $p>0$.
\begin{enumerate}
\item If $R$ is a regular ring, then it is strongly F-regular. 

\item Direct summands of F-regular rings are F-regular.

\item If $R$ is an F-rational ring, then it is normal and Cohen-Macaulay.

\item An F-rational Gorenstein ring is F-regular. 

\item The ring $R$ is F-regular if and only if it is strongly F-regular.
\end{enumerate}\label{longlist}\end{theorem}

\begin{proof} For $(1)-(4)$ see \cite[Theorem 3.1]{HHstrong} and \cite[Theorem
4.2]{HHbasec}, and for $(5)$, see \cite[Corollary 4.3]{LS}. \end{proof}

In \cite{HHchar0} Hochster and Huneke have developed a theory of tight closure
for rings essentially of finite type over fields of characteristic zero.
However, notions corresponding to F-regularity and F-rationality can also be
defined in characteristic zero without explicitly considering a closure
operation for rings of characteristic zero. A ring $R=K[X_1, \dots, X_n]/I$
over a field $K$ of characteristic zero is said to be of {\it dense F-regular
type}\/ if there exists a finitely generated $\mathbb Z$-algebra $A \subseteq
K$ and a finitely generated free $A$-algebra 
$$
R_A= A[X_1, \dots, X_n]/I_A
$$ 
such that $R \cong R_A \otimes_A K$ and, for all maximal ideals $\mu$ in a 
Zariski dense subset of $\Spec A$, the fiber rings $R_A \otimes_A A/\mu$ are 
F-regular rings of characteristic $p>0$. Similarly, $R$ is of {\it dense
F-rational type}\/ if for a dense subset of $\mu$, the fiber rings $R_A
\otimes_A A/\mu$ are F-rational. Combining some of the results from \cite{Hara,
wat-log, sm-ratsing} we have:

\begin{theorem}[Hara-Smith-Watanabe] Let $R$ be a ring finitely generated over
a field of characteristic $0$. Then $R$ has rational singularities if and only
if it is of dense F-rational type. If $R$ is $\mathbb Q$-Gorenstein, then it
has log terminal singularities if and only if it is of dense F-regular type.
\end{theorem}

When used along with Proposition \ref{effective}, the following result provides
various examples of F-rational rings which are not F-regular:

\begin{prop} Let $D$ be an ample $\mathbb Q$-divisor on ${\mathbb P}^d$, where 
${\mathbb P}^d$ is projective space over a perfect field $K$ of characteristic
$p>0$ (or of characteristic $0$), and consider the generalized section ring
$R=R({\mathbb P}^d,D)$. If $R$ is F-regular (respectively, of dense F-regular
type), then $\deg (K_{{\mathbb P}^d}+D') < 0$. \label{degree} \end{prop}

\begin{proof} If $K$ is a perfect field of characteristic $p>0$, the ring $R$ is
F-finite, i.e., the ring $R^{1/p}$ (obtained by adjoining $p\,$th roots of
elements of $R$) is a finitely generated $R$-module. By \cite{LS} $R$ is
strongly F-regular and $0^*_E = 0$, where $E$ is the injective hull of the
residue field $K$; in other words, the zero submodule of $E$ is tightly closed.
Let $\zeta \in [H^{d+1}_m(\omega_R)]_0$ denote a socle generator of $E$. Since
$0^*_E = 0$, for any nonzero element $c \in R_n$, there exists $q=p^e$ such
that the element 
$$
c F^e(\zeta) \in [H^{d+1}_m(\omega_R^{(q)})]_n
$$
is nonzero, where $F$ denotes the Frobenius morphism. Consequently 
$$
[H^{d+1}_m(\omega_R^{(q)})]_n \cong 
H^d({\mathbb P}^d,{\mathcal O}_{{\mathbb P}^d}(q(K_{{\mathbb P}^d}+D') +nD))
\neq 0,
$$
and it follows that $\deg (K_{{\mathbb P}^d}+D') < 0$. \end{proof}

\begin{example} For an algebraically closed field $K$ of characteristic $0$
and positive integers $r,\,n,$ and $m \ge 3$, consider the 
$\mathbb Q$-divisor 
$$
D = \frac{1}{r} V (X_1^n + \cdots + X_m^n) \quad \text{ on } \quad
{\mathbb P}^{m-1} = \Proj K[X_1,\dots,X_m].
$$
Note that $X_1^n + \cdots + X_m^n$ is irreducible since $m \ge 3$. Take the 
corresponding generalized section ring $R = \oplus_{n\ge 0} 
H^0({\mathbb P}^{m-1},\mathcal O_{\mathbb P^{m-1}}(nD))Z^n$
where the term $Z^n$ serves to keep track of the degree of an element of $R$.
The ring $R$ is generated over the field $K$ by $Z$, and elements of the form 
$$
\frac{\mu_n (X_1, \dots, X_m)}{X_1^n + \cdots + X_m^n}Z^r
$$
where $\mu_n(X_1, \dots, X_m)$ is a monomial of degree $n$ in the variables 
$X_1, \dots, X_m$. Consequently $R$ may be identified with the subring of the 
hypersurface
$$
S = K[Z,\,X_1,\dots,X_m]/(Z^r-(X_1^n + \cdots + X_m^n))
$$
generated by $z$, (i.e., the image of $Z$ in $S$) and monomials of degree $n$ 
in $x_1, \dots, x_m$ (the images, respectively, of $X_1,\dots,X_m$). For another
description of $R$, let $\zeta \in K$ be a primitive $n$\,th root of unity and 
consider the $K$-automorphism $\sigma$ of $S$ where $\sigma(x_i) = \zeta x_i$ 
for $1 \le i \le m$ and $\sigma(z)=z$. If $G=(\sigma) \cong 
{\mathbb Z}/n{\mathbb Z}$, then the ring of invariants of this group action, 
$S^G$, is isomorphic to $R$. It is an easy verification that $R$ has an isolated
singularity at its homogeneous maximal ideal $m_R$ since, for all $1 \le i \le
m$, the localization $R_{x_i^n}$ is a regular ring. The $\mathbb Q$-divisor $D$ 
is effective, and so the ring $R$ has rational singularities by Proposition 
\ref{effective}. If $R$ is of dense F-regular type then, by Proposition 
\ref{degree},
$$
\deg (K_{{\mathbb P}^{m-1}}+D') = -m + \left( \frac{r-1}{r} \right) n < 0.
$$
Consequently if $(r-1)n -mr \ge 0$, then $R$ is a $\mathbb Q$-Gorenstein ring, 
with an isolated rational singularity, which is not of dense F-regular type.
\label{nmr} \end{example}

\section{Cyclic covers}

For the general theory of cyclic covers, we refer the reader to the paper
of Tomari and Watanabe, \cite{tw}. The discussion below shall suffice for
our needs.

Let $R$ be an $\mathbb N$-graded or local normal domain with an ideal $I$ of 
pure height one (i.e., a {\em divisorial}\/ ideal) which has finite order $n$
when regarded as an element of the divisor class group of the ring $R$. Let
$I^{(n)}= uR$ where, in the graded case, we assume furthermore that $I$ and $u$
are homogeneous. By the {\it cyclic cover}\/ of $R$ with respect to $I$, we
mean the ring 
$$
S=R[It,I^{(2)}t^2, \dots, I^{(i)}t^i, \dots]/(ut^n-1).
$$
This ring is finitely generated as an $R$-module, in particular, the generators
of $It, \dots, I^{(n-1)}t^{n-1}$ form a generating set for $S$. If the 
characteristic of the residue field of $R$ is relatively prime to $n$, the
inclusion of $R$ in $S$ is \'etale in codimension one; more generally, 
if $P \in \Spec R$ is such that $IR_P$ is principal, say $IR_P=vR_P$, then
$u = \lambda v^n$ for a unit $\lambda \in R_P$, and so
$$ 
S_P=R_P[vt]/(\lambda (vt)^n-1), 
$$ 
which is an \'etale extension of $R_P$. Note that $S$ is a normal 
domain. We next recall a theorem of Watanabe:

\begin{theorem}\cite[Theorem 2.7]{wat-dim2} Let $R \to S$ be a finite local
homomorphism of normal local rings (or a graded homomorphism of $\mathbb
N$-graded rings) which is \'etale in codimension one. If the ring $R$ is
strongly F-regular, then the ring $S$ is also strongly F-regular. 

In particular, if $S$ is the cyclic cover of a strongly F-regular ring $R$ with
respect to a divisorial ideal of order $n$, then $S$ is strongly F-regular (and
hence Cohen-Macaulay) provided the characteristic of $R$ does not divide $n$.
\label{etale}\end{theorem}

One of the main goals of this paper is to show that under the weaker hypothesis
that the ring $R$ is F-rational, a cyclic cover $S$ need not be Cohen-Macaulay
even if $R$ is a $\mathbb Q$-Gorenstein normal domain with an isolated rational
singularity, and the finite extension $S$ is simply the canonical cover of the
ring $R$.

We focus next on the case where $R$ is an $\mathbb N$-graded normal domain, and
$I$ is a homogeneous ideal of pure height one. In this setting we show that
there is a natural $\mathbb Q$-grading on the cyclic cover $S$ which extends
the grading on $R$. The main point of our result below is that this grading
involves only non-negative rational numbers, and that the degree zero component
of $S$ is the field $[R]_0 = K$.

\begin{prop}
Let $R$ be an $\mathbb N$-graded normal domain where $[R]_0 = K$ is a field, and
$I$ be a homogeneous ideal of pure height one which has order $n$ as an element
of the divisor class group of $R$. Then $I^{(n)}=uR$ for a homogeneous element
$u \in R$, and there is a unique $\mathbb Q$-grading on the cyclic cover
$$
S=R[It,I^{(2)}t^2, \dots, I^{(i)}t^i, \dots]/(ut^n-1)
$$
which extends the grading on $R$. Under this grading, homogeneous elements of 
$S$ have nonnegative weights and $[S]_0 = [R]_0 = K$.
\label{gradedcover}\end{prop}

\begin{proof}
Since $ut^n=1$, we have $\deg t = -(\deg u)/n$ in the fraction field of $S$, 
and this gives the unique extension of the $R$-grading. Let 
$xt^i \in I^{(i)}t^i$ be a nonzero homogeneous element of $S$ where 
$1\le i\le n-1$. Then 
$$
(xt^i)^n \in I^{(ni)}t^{ni} = u^it^{ni}R = R, \qquad (*)
$$ 
and so $\deg (xt^i) \ge 0$. If $\deg (xt^i) = 0$, then $(*)$ furthermore shows 
that $(xt^i)^n = c \in [R]_0=K$, and so $x^n = cu^i$. Since $x \in I^{(i)}$, we 
have 
$$
cu^i = x x^{n-1} \in I^{(i)} I^{(in-i)} \subseteq I^{(in)}= u^iR,
$$
and so $I^{(i)} I^{(in-i)} = u^iR$. However this means that $I^{(i)}$ is an
invertible fractional ideal, and hence a projective $R$-module, 
\cite[Theorem 11.3]{mat-book}. Since $R$ and $I^{(i)}$ are graded, 
$I^{(i)}$ must be a free $R$-module, and since it has rank one, it must be
principal. However this contradicts the fact that the order of $I$ in the 
divisor class group of $R$ is $n$. \end{proof}

We usually prefer to work with this $\mathbb Q$-grading since it agrees with
the grading on $R$. However we point out that since only finitely many
denominators occur in this ${\mathbb Q}$-grading, one could multiply by a
suitable integer to obtain an ${\mathbb N}$-grading on $S$. (Specifically,
multiplying all weights by the integer $n$ will ensure an ${\mathbb
N}$-grading.) 

Our focus next is on the canonical cover of a graded $\mathbb Q$-Gorenstein
ring $R$. It is well known that if an integral domain $R$ has a canonical
module, then there is an ideal of $R$ which is isomorphic to that canonical
module $\omega$. However in the graded situation it is not always possible to
find an ideal of $R$ which is isomorphic to the graded canonical module
$\HH_K(\H^d_m(R), \, K)$ in a degree preserving manner. This obstacle may be
circumvented by choosing instead a fractional ideal of $R$ and henceforth, by a
{\em graded canonical module}\/, we shall mean a fractional ideal $\omega$ of
$R$ which is isomorphic to $\HH_K(\H^d_m(R), \, K)$ via a degree preserving
isomorphism. Since $\omega$ is a reflexive module of rank one, its symbolic
power $\omega^{(i)}$ may be defined as the double dual of the tensor product of
$i$ copies of $\omega$.

An $\mathbb N$-graded or local normal domain $R$ is said to be {\em $\mathbb
Q$-Gorenstein} if its canonical module has finite order when regarded as an
element of the divisor class group of $R$. A generalized section ring
$R(\mathbb P^d,D)$ is always $\mathbb Q$-Gorenstein: since $\Pic (\mathbb P^d)
=\mathbb Z$, there exist integers $m$ and $n$ such that $m(K_{\mathbb P^d}+D')$
and $nD$ are linearly equivalent {\em Weil}\/ divisors. By the {\it canonical
cover}\/ of a $\mathbb Q$-Gorenstein ring $R$, we mean the cyclic cover of $R$
with respect to its canonical module. For a $\mathbb Q$-Gorenstein ring $R$, we
shall use $\widetilde{R}$ to denote its canonical cover. 

We shall say that a normal $\mathbb N$-graded ring $R$ is {\em 
quasi-Gorenstein}\/ if the module $\HH_K(\H_m^d(R), \, K)$ and the ring $R$ are 
isomorphic as underlying $R$-modules. In this case, it follows that
$$
\omega_R = (\H_m^d(R))^* \cong R(a)
$$ 
where $a$ is the $a$-invariant of the ring $R$. A quasi-Gorenstein \CM ring 
is Gorenstein.

\begin{prop}
Let $R$ be a graded normal $\mathbb Q$-Gorenstein ring. If the canonical module 
$\omega$ of $R$ has order $n$ in the divisor class group of $R$ and 
$\omega^{(n)} = uR$, consider the graded canonical cover 
$$
\widetilde{R}=R[\omega t,\omega^{(2)}t^2,\dots,\omega^{(n-1)} t^{n-1}, \dots]
/ (ut^n-1)
$$
as above, i.e., with $\deg t= -(\deg u)/n$. Then the ring $\widetilde{R}$ 
is quasi-Gorenstein with $a$-invariant $a(\widetilde{R}) = -(\deg u)/n$. 
\label{a-inv} \end{prop}

\begin{proof}
If $k = (\deg u)/n$, we have a graded isomorphism
$$
\widetilde{R} \cong R \oplus \omega(k) \oplus \omega^{(2)}(2k) \oplus \dots 
\oplus \omega^{(n-1)}(nk-k). 
$$
If $d=\dim R$ and $0 \le i \le n-1$, there are graded isomorphisms
\begin{align*}
\HH_K(\H_m^d(\omega^{(i)}), \, K) & \cong
 \HH_K(\H_m^d(\omega) \, \HT_R \, \omega^{(i-1)}, \, K) \\
& \cong \HH_R(\omega^{(i-1)}, \, \HH_K(\H_m^d(\omega), \, K)) \\
& \cong \HH_R(\omega^{(i-1)}, \, R) \ \cong \ \omega^{(1-i)}.
\end{align*}
Consequently 
\begin{align*}
\omega_{\widetilde{R}} \ & \cong \ \HH_K(\H_m^d(\widetilde{R}), \, K) 
\cong \ \bigoplus_{i=0}^{n-1} \HH(\H_m^d(\omega^{(i)})(ik), \, K) \\
& \cong \ \bigoplus_{i=0}^{n-1} \omega^{(1-i)}(-ik) \ \cong \ \widetilde{R}(-k),
\end{align*}
and so $\widetilde{R}$ is quasi-Gorenstein with $a$-invariant 
$a(\widetilde{R}) = -k = -(\deg u)/n$. 
\end{proof}

We record an example we shall use later:

\begin{example} Let $K$ be an algebraically closed field of characteristic $0$, 
and take the $\mathbb Q$-divisor
$$
D=\frac{1}{3}V(y_0)+\frac{1}{3}V(z_0)+\frac{1}{3}V(y_0+z_0) \quad \text{ on } 
\quad {\mathbb P}^1 = \Proj K[y_0,\,z_0]
$$
where, for example, $V(y_0)$ denotes the point of $\mathbb P^1$ defined by the 
vanishing of $y_0$. Setting $A=\oplus_{n\ge 0} 
H^0({\mathbb P^1},\mathcal O_{\mathbb P^1}(nD))x^n$, the generators of $A$ are
$$
\frac{y_0^3x^3}{y_0z_0(y_0+z_0)}, 
\ \frac{y_0^2z_0x^3}{y_0z_0(y_0+z_0)},
\ \frac{y_0z_0^2x^3}{y_0z_0(y_0+z_0)},
\ \frac{z_0^3x^3}{y_0z_0(y_0+z_0)},
\quad \text{and} \quad x,
$$
and so we may identify $A$ with the subring of the hypersurface
$$
H = K[X,\,Y,\,Z]/(X^3-YZ(Y+Z))
$$
generated by the elements $y^3,\,y^2z,\,yz^2,\,z^3$ and $x$. A graded canonical 
module for $A$ is 
$$
\omega_A = \frac{1}{x^2}\left( y^3,\,y^2z \right)A,
$$ 
and its symbolic powers are
$$
\omega_A^{(2)} = \frac{1}{x^4}\left(y^6,\,y^5z,\,y^4z^2\right)A
\qquad \text{ and } \qquad \omega_A^{(3)} = \frac{1}{x^6}\left(y^6\right)A.
$$
The isomorphism defined by
$$
\frac{y^3t}{x^2} \mapsto y, \ \frac{y^2zt}{x^2} \mapsto z, 
\ \frac{y^6t^3}{x^6} \mapsto 1
$$
may be used to identify $\widetilde{A}$ with the hypersurface $H$. \label{333}
\end{example}

We next examine when the Segre product of two graded $\mathbb Q$-Gorenstein
rings is a $\mathbb Q$-Gorenstein ring.

\begin{prop} Let $A$ and $B$ be normal $\mathbb N$-graded $\mathbb
Q$-Gorenstein rings over an algebraically closed field $[A]_0 = [B]_0 = K$. If
the graded canonical modules $\omega_A$ and $\omega_B$ have orders $m$ and $n$
respectively as divisor class group elements, choose homogeneous elements $u$
and $v$ in the fraction fields of $A$ and $B$ such that $\omega_A^{(m)} \cong
uA$ and $\omega_B^{(n)} \cong vB$ as graded fractional ideals. If 
$n\deg u = m\deg v$, then the Segre product $R = A \# B$ is also a 
$\mathbb Q$-Gorenstein ring. If, furthermore, $m$ and $n$ are relatively prime 
integers, then the Segre product of the graded canonical covers of $A$ and $B$ 
is isomorphic to the canonical cover of $R=A \# B$, i.e.,
$$
\widetilde{A \# B} \ \cong \ \widetilde{A} \# \widetilde{B}.
$$ 
\label{cov-seg}
\end{prop} 

\begin{proof}
It is proved in \cite{GW} that the canonical module $\omega_R$ is the Segre
product of the graded canonical modules, $\omega_A \#\omega_B$, see also \S\,2. 
We have
$$
\omega_R^{(mn)} \cong \omega_A^{(mn)} \# \omega_B^{(mn)} 
\cong u^n A \# v^m B
$$
and this is isomorphic to $u^n v^m (A \# B)$ whenever $n\deg u = m\deg v$.

In the case that $m$ and $n$ are relatively prime, $\omega_R$ has order $mn$ as 
a divisor class group element. If
$$
k = (\deg u)/m = (\deg v)/n, 
$$ 
then the graded canonical cover of $R$ is 
$$
\widetilde{R} \ \cong \ \bigoplus_{r=0}^{mn-1} \omega_R^{(r)}(kr) \ \cong \ 
\bigoplus_{r=0}^{mn-1} \omega_A^{(r)}(kr) \ \# \ \omega_B^{(r)}(kr). 
$$
Since $m$ and $n$ are relatively prime integers, given arbitrary $i,\,j \in 
{\mathbb N}$, there exists a unique integer $r$ with $0 \le r \le mn-1$ such 
that $r \equiv i \mod m$ and $r \equiv j \mod n$. Consequently
$$
\widetilde{R} \ \cong \ \bigoplus_{i=0}^{m-1} \ \bigoplus_{j=0}^{n-1} 
\omega_A^{(i)}(ki) \ \# \ \omega_B^{(j)}(kj) \ 
\cong \ \widetilde{A} \# \widetilde{B}.
$$
\label{q-gor-segre}\end{proof}

By Watanabe's result, Theorem \ref{etale}, the canonical cover of a $\mathbb
Q$-Gorenstein strongly F-regular ring $R$ is also strongly F-regular, provided
the characteristic of $R$ does not divide the order of $[\omega_R] \in \Cl(R)$.
The corresponding result for log terminal singularities was proved earlier by
Kawamata:

\begin{theorem}\cite[Proposition 1.7]{Ka} Let $R$ be a normal $\mathbb
Q$-Gorenstein ring of characteristic zero with a canonical cover
$\widetilde{R}$. Then $R$ has log terminal singularities if and only if
$\widetilde{R}$ has log terminal singularities. \label{etale2}\end{theorem}

\subsection*{Infinite covers}

When the canonical ideal $\omega_R$ of a \CM normal ring $(R,m)$ is not
necessarily of finite order in the divisor class group $\Cl(R)$, the {\it
anti-canonical cover}\/ $S$ may be constructed by taking an ideal $I$ of pure 
height one which is an inverse for $\omega_R$ in $\Cl(R)$, and forming the 
symbolic Rees ring 
$$
S=\oplus_{i \ge 0}I^{(i)}.
$$
The interest in the anti-canonical cover arises from the fact that if $S$ as
above is Noetherian and Cohen-Macaulay, then it is Gorenstein. Symbolic Rees
rings in general need not be Noetherian, but there is an interesting theorem of
Watanabe that applies when the symbolic Rees ring is Noetherian:

\begin{theorem}\cite[Theorem 0.1]{wat-antic} Let $(R,m)$ be a strongly
F-regular ring, and $I$ an ideal of pure height one which is the inverse of the
canonical module $\omega_R$ in $\Cl(R)$. If the anti-canonical cover
$S=\oplus_{i \ge 0}I^{(i)}$ is Noetherian, then it is strongly F-regular.
\label{antican} \end{theorem}

For an arbitrary divisorial ideal $I$ of a strongly F-regular ring $R$,
Watanabe raised the question whether the symbolic Rees algebra $\sym(I) =
\oplus_{n \ge 0}I^{(n)}$ is Cohen-Macaulay whenever it is Noetherian. In
\cite{singh-rees} we developed the notion of multi-symbolic Rees algebras and
used it to show that $\sym(I)$ is indeed \CM whenever a certain auxiliary ring
is finitely generated over $R$. More precisely, we established:

\begin{theorem}\cite[Theorem 5.1]{singh-rees}
Let $(R,m)$ be a strongly F-regular ring with canonical ideal $\omega_R$.
Given a divisorial ideal $I$, choose a divisorial ideal $J$ such that 
$[I]+[J]+[\omega_R]=0$ in the divisor class group $\Cl(R)$. If the
multi-symbolic Rees algebra 
$$
\sym(I,J) = \bigoplus_{n,m \ge 0} \HR ( \HR ( I^nJ^m,R),R)
$$ 
is finitely generated over $R$, then $\sym(I)$ is strongly F-Regular and, in
particular, is Cohen-Macaulay. 
\end{theorem}

Note that $\HR (\HR (I^nJ^m,R),R)$ above is simply the reflexive hull of
$I^nJ^m$. The hypothesis that $R$ is strongly F-regular is indeed used in an
essential way: Watanabe has constructed an example of an F-rational ring $R$
with a divisorial ideal $I$ such that the symbolic Rees algebra $\sym(I)$ is
not Cohen-Macaulay, \cite[Example 4.4]{wat-antic}. Our result Theorem \ref{ncm}
may be viewed as a strengthening of this in the sense that we obtain such an
example with the additional restriction that the divisorial ideal $I$ has
finite order in $\Cl(R)$.

\section{Griffith's examples from a tight closure viewpoint} \label{phil}

Let $f: X \to Y$ be a finite morphism of smooth projective varieties over an
algebraically closed field $K$ of characteristic zero, and $\mathcal L$ be an 
ample line bundle on $Y$ such that the homogeneous coordinate ring 
$\oplus_{n \ge 0} H^0(Y, {\mathcal L}^{\otimes n})$ is Cohen-Macaulay. In 
\cite{griffith} Griffith examined the issue whether the arithmetically \CM 
property ascended to $X$, i.e., whether $\oplus_{n \ge 0} H^0(X, 
f^* {\mathcal L}^{\otimes n})$ is Cohen-Macaulay, under the hypothesis that the 
homomorphism
$$
\oplus_{n \ge 0} H^0(Y,{\mathcal L}^{\otimes n}) \to 
\oplus_{n \ge 0} H^0(X, f^* {\mathcal L}^{\otimes n}) 
$$
is \'etale away from the irrelevant ideal. Griffith constructed examples of
morphisms $f: X \to Y$ where the arithmetically Cohen-Macaulay property ascends
but, if $r \ge 1$, it does not ascend for the induced morphism $f: X \times
{\mathbb P}^r \to Y \times {\mathbb P}^r$.

We recall the examples from \cite{griffith}, and discuss them from the point of
view of log terminal singularities. For $d \ge 4$, let $S$ be the hypersurface
$$
S= K[Z,\,X_1,\dots,X_{d-1}]/(Z^d-(X_1^d+\cdots+X_{d-1}^d))
$$
over an algebraically closed field $K$ of characteristic zero. Let $\zeta$ be a
primitive $d$\,th root of unity and consider the $K$-automorphism $\sigma$ of 
$S$ where $\sigma(x_i) = \zeta x_i$ for $1 \le i \le d-1$ and $\sigma(z)=z$. 
Let $G=(\sigma) \cong {\mathbb Z}/d{\mathbb Z}$. Then the ring of invariants of
this group action $R=S^G$ is generated over $K$ by $z$ and the monomials of
degree $d$ in the elements $x_1,\dots,x_{d-1}$. Let $A=K[Y_0,\dots,Y_r]$ be a
polynomial ring over $K$ where $r \ge 1$ and set $X=\Proj S$, $Y=\Proj R$, and 
${\mathbb P}^r = \Proj A$. In his paper, Griffith established the 
following results:

\begin{prop}[\cite{griffith} \S\,4] With the above notation, 
\begin{enumerate}
\item The rings $R$ and $S$ are \CM and the extension $R \to S$ is generically
Galois, with Galois group ${\mathbb Z}/d{\mathbb Z}$. The extension of Segre
products $R\# A \to S\# A$ is also generically Galois, with Galois group 
${\mathbb Z}/d{\mathbb Z}$.

\item The graded $K$-algebra $R$ satisfies $a(R)<0$ whereas $a(S)=0$.
Consequently the ring $R\# A$ is Cohen-Macaulay, but $S\# A$ is not. More
specifically, $S\# A$ satisfies the Serre condition $S_{d-1}$ but not the
condition $S_d$.

\item The extension $R \to S$ is \'etale away from the irrelevant maximal
ideal, as is the extension $R\# A \to S\# A$.
\end{enumerate}\end{prop}

\begin{remark} The ring $R$ above is the case of Example \ref{nmr} with $r=n=d$ 
and $m=d-1$, i.e., it is isomorphic to the generalized section ring of the 
$\mathbb Q$-divisor 
$$
D = \frac{1}{d} V (X_1^d + \cdots + X_{d-1}^d) \quad \text{ on } \quad
\Proj K[X_1, \dots, X_{d-1}].
$$
As a graded canonical module for $R$, we may take
$$
\omega_R = \frac{1}{z^{d-1}} \big(x_1^d,\,x_1^{d-1}x_2,\dots,x_1^{d-1}x_{d-1})R.
$$
This has order $d$ as an element of $\Cl(R)$, and $\omega_R^{(d)} =
(x_1^{d(d-1)}/z^{d(d-1)})R$. Furthermore, it turns out that the canonical cover
of $R$ is isomorphic to the hypersurface $S$. Note that $\deg (K_{\mathbb
P^{d-2}} + D') = 0$ and so $R$ is an example of a $\mathbb Q$-Gorenstein ring
with an isolated rational singularity, which is not a log terminal singularity. 
Such a ring always yields an example in which the arithmetically Cohen-Macaulay
property fails to ascend, as we record in the proposition below: \end{remark}

\begin{prop} Let $R$ be an $\mathbb N$-graded $\mathbb Q$-Gorenstein ring of
characteristic zero with an isolated rational singularity which is not a 
log terminal singularity. Then at least one of the following two conditions 
holds:
\begin{enumerate}
\item The canonical cover $\widetilde{R}$ of the ring $R$ is not
Cohen-Macaulay, in other words, the arithmetically \CM property does not ascend
under the morphism $\Proj \widetilde{R} \to \Proj R$.

\item The canonical cover satisfies $a(\widetilde{R}) \ge 0$ and so, if $r \ge
1$, the arithmetically \CM property does not ascend under the morphism 
$$
(\Proj \widetilde{R}) \times {\mathbb P}^r \to (\Proj R) \times {\mathbb P}^r.
$$
\end{enumerate}\end{prop}

\begin{proof} By Kawamata's result, recorded earlier as Theorem \ref{etale2},
the canonical cover $\widetilde{R}$ does not have log terminal singularities.
But then, since $\widetilde{R}$ is quasi-Gorenstein, it cannot have rational
singularities. By our assumption $R$ has an isolated singularity, and therefore
the canonical cover $\widetilde{R}$ has an isolated singularity as well. By
Watanabe's result, Theorem \ref{ratsing-ainv}, the failure of $\widetilde{R}$
to have rational singularities can be attributed to one of the following two
reasons:

$(1)$ \quad The ring $\widetilde{R}$ is not Cohen-Macaulay, and so the
arithmetically Cohen-Macaulay property does not ascend under the morphism
$\Proj \widetilde{R} \to \Proj R$.

$(2)$ \quad The ring $\widetilde{R}$ has a nonnegative $a$-invariant. In this
case, consider the polynomial ring $A=K[Y_0,\dots,Y_r]$ where $r \ge 1$. Then
$R\#A$ is a Cohen-Macaulay ring whereas, by the K\"unneth formula for local
cohomology, $\widetilde{R}\#A$ is not Cohen-Macaulay. Hence the arithmetically
Cohen-Macaulay property does not ascend under the morphism $(\Proj
\widetilde{R}) \times {\mathbb P}^r \to (\Proj R) \times {\mathbb P}^r$ for $r
\ge 1$. \end{proof}

The canonical cover $\widetilde{R}$ of the ring $R$ constructed by Griffith is 
isomorphic to the hypersurface 
$$
S=K[X_1,\dots,X_{d-1},Z]/(Z^d-(X_1^d+\cdots+X_{d-1}^d))
$$
which, while being Cohen-Macaulay, does not have rational singularities since
$a(S) = 0$, i.e., due to reason $(2)$ above. The question then arises: does
there exist a ring $R$, with an isolated rational singularity, whose canonical
cover $\widetilde{R}$ does not have rational singularities for reason $(1)$,
i.e, because $\widetilde{R}$ is not Cohen-Macaulay? The answer is {\em yes}, as
we see next.

\section{Non Cohen-Macaulay canonical covers of rings with rational
singularities}

We show that the canonical cover of a ring with rational singularities need not
be Cohen-Macaulay:

\begin{theorem} For all $d \ge 3$, there exists a graded $\mathbb Q$-Gorenstein
ring $R$ of dimension $d$, with an isolated rational singularity, whose
canonical cover $\widetilde{R}$ has depth $2$. In particular, $\widetilde{R}$
is not a Cohen-Macaulay ring. \label{ncm} \end{theorem}

\begin{proof} Let $K$ be an algebraically closed field of characteristic $0$. 
The ring $R$ in the statement of the theorem will be constructed as a Segre
product $A\# B$. For the ring $A$, as in Example \ref{333}, take the subring of 
$$
\widetilde{A} = K[X,\,Y,\,Z]/(X^3-YZ(Y+Z))
$$
generated by the elements $y^3,\,y^2z,\,yz^2,\,z^3$ and $x$. 

{\bf Case $d=3:$} \quad On ${\mathbb P}^1 = \Proj K[v_0,\,w_0]$, take the 
$\mathbb Q$-divisor
$$
E=\frac{1}{2}V(v_0)+\frac{1}{2}V(w_0)+\frac{1}{2}V(v_0-w_0)
+\frac{1}{2}V(v_0+w_0)
$$ 
and set $B = \oplus_{n \ge 0} 
H^0({\mathbb P^1},\mathcal O_{\mathbb P^1}(nE))u^n$. 
The generators of $B$ are
\begin{align*}
& \frac{v_0^4 u^2}{v_0w_0(v_0^2-w_0^2)}, 
\ \frac{v_0^3 w_0 u^2}{v_0w_0(v_0^2-w_0^2)}, 
\ \frac{v_0^2 w_0^2 u^2}{v_0w_0(v_0^2-w_0^2)}, 
\ \frac{v_0 w_0^3 u^2}{v_0w_0(v_0^2-w_0^2)}, \\
& \frac{w_0^4 u^2}{v_0w_0(v_0^2-w_0)}, \quad \text{ and } \quad u,
\end{align*} and we identify $B$ with the subring of the hypersurface
$$
K[U,\,V,\,W]/(U^2-VW(V^2-W^2))
$$
generated by the elements $v^4,\,v^3w,\,v^2w^2,\,vw^3,\,w^4$, which are 
assigned weight $2$, and $u$ which has weight $1$. A graded canonical module 
for $B$ is 
$$
\omega_B = \frac{1}{u}\left(v^4,\,v^3w,\,v^2w^2\right)B.
$$ 
This has order $2$ as a divisor class group element, and 
$\displaystyle{\omega_B^{(2)} = \frac{v^4}{u^2}B}$. The canonical cover 
$\widetilde{B}$ is isomorphic to the subring of the hypersurface generated by 
$v^2,\,vw,\,w^2$ and $u$.

{\bf Case $d \ge 4:$} \quad In the situation of Example \ref{nmr}, take
$r=2$, $m=d-1$, and $n=2m$, i.e., the $\mathbb Q$-divisor
$$
E = \frac{1}{2} V (X_1^{2m} + \cdots + X_m^{2m}) \quad \text{ on } \quad
{\mathbb P}^{m-1} = \Proj K[X_1,\dots,X_m],
$$
and the generalized section ring $B = \oplus_{n\ge 0} H^0({\mathbb P}^{m-1}, 
\mathcal O_{\mathbb P^{m-1}}(nE))Z^n$. Note that $B$ is not of dense 
F-regular type, since $\deg (K_{{\mathbb P}^{m-1}}+E') = 0$. We identify $B$ 
with the subring of the hypersurface 
$$
K[Z,\,X_1,\dots,X_m]/(Z^2-(X_1^{2m} + \cdots + X_m^{2m})), 
$$
generated by $z$ and monomials of degree $2m$ in $x_1,\dots,x_m$. It is easily 
checked (using the corresponding $\mathbb Q$-divisors) that as a graded 
canonical module for $B$, we may take the fractional ideal
$$
\omega_B = \frac{1}{z} \big( x_1^m\mu_m \ : \ \mu_m
\text{ runs through degree $m$ monomials in } x_1,\dots,x_m \big)B.
$$
Its second symbolic power is $\omega_B^{(2)} = (x_1^{2m}/z^2)B$.

Note that the order of $\omega_A$ in $\Cl(A)$ is $3$, and (in either case) that
$\omega_B$ has order $2$ as an element of $\Cl(B)$. Furthermore,
$a(\widetilde{A}) = a(\widetilde{B}) =0$, so by Proposition \ref{cov-seg}, the
Segre product $R = A \# B$ is $\mathbb Q$-Gorenstein, and its canonical cover
$\widetilde{R}$ is isomorphic to $\widetilde{A} \# \widetilde{B}$. By
Proposition \ref{effective}, $A$ and $B$ are rings with rational singularities
over an algebraically closed field $K$, it follows that $A \otimes_K B$ has
rational singularities. By Boutot's Theorem, \cite{boutot}, its direct summand
$R = A \# B$ has rational singularities as well. The verification that $R$ has
an isolated singularity at $m_R$ is rather routine, and is left to the reader.

It remains to check that the canonical cover of $R$, i.e., $\widetilde{R} \cong 
\widetilde{A} \# \widetilde{B}$, has depth $2$. Since $\dim \widetilde{A} = 2$, 
by the K\"unneth formula for local cohomology we have,
$$
H_{m_{\widetilde{R}}}^2(\widetilde{R}) \cong \left(
\widetilde{A} \# H_{m_{\widetilde{B}}}^2(\widetilde{B}) \right) \oplus \left( 
H_{m_{\widetilde{A}}}^2(\widetilde{A}) \# \widetilde{B} \right)
$$
which is nonzero since $a(\widetilde{A})=0$.
\end{proof}

\bibliographystyle{amsalpha}

\end{document}